\documentclass{amsart}
\usepackage{amsfonts}
\usepackage{amsmath}
\usepackage{graphicx}
\usepackage{url}
\usepackage{bm}
\usepackage{hyperref}

\newtheorem{theorem}{Theorem}[section]
\newtheorem{lemma}[theorem]{Lemma}

\theoremstyle{definition}
\newtheorem{definition}[theorem]{Definition}

\newtheorem{corollary}[theorem]{Corollary}

\numberwithin{equation}{section}

\begin{document}

\title[The Weighted SCFE Problem in the Line]{The Weighted Sitting Closer to Friends than Enemies Problem in the Line}

\author{Julio Aracena}\address{CI$^2$MA \and Departamento de Ingenier\'ia Matem\'atica, Facultad de Ciencias F\'isicas y Matem\'aticas, Universidad de Concepci\'on, Chile.}\email{jaracena@ing-mat.udec.cl} \thanks{Julio Aracena was partially supported by ANID-Chile through the project {\sc Centro de Modelamiento Matem\'atico} (AFB170001) of the PIA Program:``Concurso Apoyo a Centros Cient\'ificos y Tecnol\'ogicos de Excelencia con Financiamiento Basal'', and by ANID-Chile through Fondecyt project 1151265.}

\author{Christopher Thraves Caro}\address{Departamento de Ingenier\'ia Matem\'atica, Facultad de Ciencias F\'isicas y Matem\'aticas, Universidad de Concepci\'on, Chile.}\email{cthraves@ing-mat.udec.cl}


\date{}


\begin{abstract}
The weighted \emph{Sitting Closer to Friends than Enemies} (SCFE) problem is to find an injection of the vertex set of a given weighted
graph into a given metric space so that, for every pair of incident edges with different weight, the end vertices of the heavier edge are closer 
than the end vertices of the lighter edge.
The \emph{Seriation} problem is to find a simultaneous reordering of the rows and columns of a symmetric matrix such that the entries are monotone nondecreasing in rows and
columns when moving towards the diagonal. If such a reordering exists, it is called a \emph{Robinson} ordering.   
In this work, we establish a connection between the SCFE problem and the Seriation problem.
We show that if the \emph{extended adjacency matrix} of a given weighted graph $G$ has no Robinson ordering then $G$ has no injection in $\mathbb{R}$ that solves the SCFE problem. 
On the other hand, if the extended adjacency matrix of $G$ has a Robinson ordering, 
we construct a polyhedron that is not empty if and only if there is an injection of the vertex set of $G$ in $\mathbb{R}$ that solves the SCFE problem.
As a consequence of these results, 
we conclude that 
deciding
the existence of (and constructing) such an injection in $\mathbb{R}$ for a given 
\emph{complete}
weighted graph can be done in polynomial time. 
On the other hand, we show that deciding if an \emph{incomplete} 
weighted graph has such an injection in $\mathbb{R}$
is NP-Complete.
\end{abstract}

\maketitle

\section{Introduction}\label{sec:intro}
Consider a data set. The task is to construct a graphic 
representation of
the data set so that similarities between data points are graphically 
expressed.
To complete this task, the only information
available is a \emph{similarity matrix} of the data set, i.e., a 
symmetric, square matrix 
whose entry $ij$ contains a similarity measure between data points 
$i$ and $j$ (the larger the value the more similar the data points are).
Hence, the task is to draw all data
points in a \emph{paper} so that for
every 
three data points $i$, $j$, and $k$, if 
$i$ 
is at least as similar to $j$ than $k$,
then $i$ should be placed closer in the drawing to $j$ than $k$.
In colloquial words, for each data point $j$, the farther the other 
data points are, the less similar they are to $j$.
 
A slightly simpler version of this problem, introduced in 
\cite{kermarrec2011can}, has been studied under the name of the Sitting 
Closer to Friends than Enemies (SCFE) problem. 
The SCFE problem uses signed graphs as an input. Therefore, the 
similarity matrix 
has entries $1$ and $-1$, representing similarity and dissimilarity, or 
friendship and enmity between the data points, from where 
the problem obtains its name. The SCFE problem has been studied in the 
real line \cite{kermarrec2011can,cygan2015sitting,pardo2015embedding} 
and in the circumference 
\cite{benitez2018sitting} (which means that the \emph{paper} is
the real line or the circumference). In both cases, the real line and
the
circumference, it has been shown 
that deciding the existence of such an injection for a given signed 
graph
is NP-Complete.
Nevertheless, in both cases again, when the problem is restricted to 
complete signed graphs there exists a characterization of the families 
of complete signed graphs 
that admit a solution for the SCFE problem and it can be decided in 
polynomial time \cite{kermarrec2011can,benitez2018sitting}. 
Therefore, a natural next step is to consider the case when 
similarities range in an extended set of values. 

The SCFE problem in the line seems to be 
closely related to the \emph{Seriation} problem.
Liiv in \cite{liiv2010seriation} defines the Seriation 
problem as ``an exploratory data analysis 
technique to reorder objects into a 
sequence along a one-dimensional continuum so that it best reveals 
regularity and pattering among the whole series''. 
Seriation has applications in archaeology \cite{petrie1899sequences}, 
data visualization \cite{brusco2006branch}, exploratory analysis 
\cite{hubert2001combinatorial},
bioinformatics \cite{tien2008methods}, and machine learning 
\cite{ding2004linearized}, among others. Liiv in 
\cite{liiv2010seriation}
presents an interesting survey on 
seriation, matrix reordering and its applications.
The first important contribution of this document is to 
show that the SCFE and the Seriation problems are different. Indeed, 
we show that seriation is a necessary condition to 
solve the SCFE problem, but 
it is not sufficient. 

To continue with our exposition, in Section \ref{sec:def} we introduce the notation and
definitions used along the document. 
The rest of the document is organized as follows. 
In Section \ref{sec:relwork}, we present the state of the art and 
contextualize our contributions.
In Section \ref{sec:realline}, we present the characterization 
of weighted graphs with an injection in $\mathbb{R}$ that satisfies the
restrictions of the SCFE problem. Furthermore, we present the 
results related with complete weighted graphs. 
In Section \ref{sec:incompletegraphs}, 
we present the results regarding incomplete weighted graphs.
We conclude in Section \ref{sec:conclusions} with some final remarks and
future work.  
\section{Notation and Definitions}\label{sec:def}
We use standard notation. A graph is denoted by $G=(V,E)$. We 
consider only undirected graphs, without parallel edges and loopless. 
The set of vertices of $G$ is $V$ and the set of edges is $E$, a set of 
$2$-elements subsets of $V$. We use $n$ and $m$ to denote $|V|$ and 
$|E|$, respectively. 
Two distinct vertices $u$ and $v$ in $V$ are said to be 
\emph{neighbors} if 
$\{u,v\} \in E$. In that case, 
we say that they  are connected by an edge
which is denoted by $\{u,v\}$. 
A graph is said to be \emph{complete} if 
every pair of distinct vertices is connected by an edge, otherwise, we say that 
it is \emph{incomplete}. 

In this document, we work with weighted graphs. We denote by 
$w:E\rightarrow \mathbb{R^+}$ 
a positive real valued function that assigns $w(\{u,v\})$, a positive 
real value, to the edge $\{u,v\}$ in $E$. 
For our purposes, we consider that $w$ is a similarity measure, i.e., 
for any $\{u,v\}\in E$ the value $w(\{u,v\})$ measures how similar $u$ 
and $v$ are.
We consider that the larger the similarity measure is, the more similar 
the vertices are. 
It is worth mentioning that the fact that the weights are positive is 
just a choice made for simplicity. Actually, if we have negative weights, one can 
translate all weights by a constant and obtain only positive weights. 
All our results are still valid if we remove this assumption.

Let $(\mathcal{M}, d)$ be a metric space. 
A \emph{drawing} of a graph $G=(V,E)$ into  $\mathcal{M}$ is an 
injection $D:V\rightarrow \mathcal{M}$. 
We define a certain type of drawings that  capture the requirements of 
the SCFE problem.   
\begin{definition}\label{def:validdrawing}
Let $G=(V,E)$ be a graph, and $w:E\rightarrow \mathbb{R^+}$ be a 
positive function on $E$. 
Let $(\mathcal{M}, d)$ be a metric space.
We say that a drawing $D$ of $G$ into $\mathcal{M}$ is \emph{valid 
distance} if,
for all pair $\{t,u\}$, $\{t,v\}$ of incident edges in $E$ such that  
$w(\{t,u\}) > w(\{t,v\})$,
\[
 d(D(t),D(u)) < d(D(t),D(v)).
\]
\end{definition}

In colloquial words, a drawing is valid distance, or simply 
\emph{valid}, when 
it places vertices $t$ and $u$ strictly closer than $t$
and $v$ in $\mathcal{M}$ 
whenever $t$ and $u$ have a strictly larger similarity measure 
than $t$ and $v$.  
Now, the weighted SCFE problem in its most general presentation is 
defined as follows. 
\begin{definition}
Given a weighted graph $G$ and a metric space $\mathcal{M}$, the 
weighted
SCFE problem in $\mathcal{M}$ is to decide whether $G$ has a valid 
drawing in 
$\mathcal{M}$, and, 
in case of a positive answer, find one.
\end{definition}

In this document, we focus our attention on the case when the metric 
space is the real line, i.e.,
we consider the metric space to be the set of real values $\mathbb{R}$ 
with the Euclidean distance.  

Since we present a matrix oriented analysis, we introduce the next two 
matrix related definitions. 
Given a matrix $A$, the entry in the $i$-\emph{th} row and $j$-\emph{th}
column of $A$ is denoted by $A_{ij}$.  
For every weighted graph $G=(V,E)$ and an ordering $\pi$ of the vertex set $V$, 
we define $A^\pi(G)$, the \emph{similarity matrix} of $G$ ordered according to $\pi$, as follows. Let $\pi_i$ be the $i$-\emph{th}
element of $V$ according to $\pi$, then: 
\[
A^\pi(G)_{ij} = 
\begin{cases}
* & \text{ if } i \neq j \mbox{ and } \{\pi_i,\pi_j\} \notin E, \\
w(\{\pi_i,\pi_j\}) & \mbox{ if } i \neq j \mbox{ and } \{\pi_i,\pi_j\} \in E, \\
\max_{e \in E}w(e) & \mbox{ if } i = j.
\end{cases}
\]
This matrix is also known as 
the extended weighted \emph{adjacency matrix} of $G$.
The $i$-\emph{th} row (and $i$-\emph{th} column) contains the 
similarities between vertex $\pi_i$ and the rest of the vertices of $G$.
We may use only $A^\pi$ when the graph $G$ is contextually clear or only $A$ when $G$ and $\pi$ are
contextually clear. 
Note that any similarity matrix of any weighted graph is symmetric since
$w$ is symmetric.
A similarity matrix of a complete weighted graph does not have 
entries with the symbol $*$. 
In that case, we say that a similarity matrix is \emph{complete}, 
otherwise we say that it is \emph{incomplete}. 

W. S. Robinson in \cite{robinson_1951} introduced Robinsonian matrices. 
A complete similarity matrix $A$ is said to be \emph{Robinson}
if its entries are monotone nondecreasing in rows and 
columns when moving towards the diagonal, i.e., if for all
integers $1\leq i < j \leq n$, 
\[
A_{ij} \leq \min\{A_{ij-1},A_{i+1j}\}.
\]
Equivalently, a complete similarity matrix $A$ is Robinson if for all 
integers $1\leq i < l \leq n$ and $j,k \in [i,l]$:
\[
A_{il} \leq \min\{A_{ij},A_{kl}\}.
\]
On the other hand, a complete similarity matrix $A$ is \emph{Robinsonian}  
if its rows and columns 
can be reordered simultaneously by a permutation $\pi$ such that $A^\pi$ is Robinson. 

The Robinson matrix definition can be naturally extended to incomplete 
matrices.  In that case, a similarity matrix is \emph{Robinson} 
if its specified entries are monotone nondecreasing in rows and 
columns when moving towards the diagonal, 
i.e., 
an incomplete similarity matrix $A$ is Robinson
if for all 
integers $1\leq i < l \leq n$ and $j,k \in [i,l]$, such that  
$A_{il} \neq *$, $A_{ij} \neq *$ and $A_{kl} \neq *$,
\[
A_{il} \leq \min\{A_{ij},A_{kl}\}.
\]

Again, we say that a similarity matrix is \emph{Robinsonian}  if
its rows and columns can be simultaneously reordered in such a way that 
it passes to be Robinson.
Finally, we say that an ordering $\pi$ of the vertex set $V$ of a weighted graph $G=(V,E)$ 
is Robinson if $A^\pi(G)$ is Robinson. 

\section{Context, Related Work, and Our 
Contributions}\label{sec:relwork}
Robinsonian matrices were defined by W. S. Robinson in  
\cite{robinson_1951} in a study on how to order chronologically 
archaeological deposits. The \emph{Seriation} problem introduced in the 
same work  
is to decide whether the similarity matrix of a data set 
is Robinsonian and
reorder it as a Robinson matrix if possible. 
Recognition of complete Robinsonian matrices has been studied by several
authors. Mirkin et al. in \cite{graphsandgenes84} presented an $O(n^4)$ 
recognition algorithm, where $n\times n$ is the 
size of the matrix. On the other hand, using divide and conquer 
techniques, Chepoi et al. in \cite{chepoi1997recognition} introduced an 
$O(n^3)$ recognition algorithm. Later,
Pr\'ea and Fortin in \cite{prea2014optimal} provided an $O(n^2)$ optimal
recognition algorithm for complete Robinsonian matrices using PQ trees. 

Using the relationship between Robinsonian matrices and unit interval 
graphs presented in \cite{roberts69}, Monique Laurent and Matteo 
Seminaroti  in \cite{laurent2017lex}
introduced a recognition algorithm for Robinsonian matrices that uses 
Lex-BFS, whose time complexity is 
$O(L(m+n))$, where $m$ is the number of nonzero entries in the 
matrix,
and $L$ is the number of different 
values in the matrix. 
Later in \cite{laurent2017similarity}, the same authors presented a 
recognition algorithm with time complexity $O(n^2+nm\log n )$
that uses similarity-first search, an extension of breadth-first search 
for weighted graphs. 
Again, using the relationship between Robinsonian matrices and unit 
interval graphs, Laurent et al. in \cite{laurent2017structural}
gave a characterization of Robinsonian matrices via forbidden patterns.

The Seriation problem also has been studied as an optimization problem.
Given an  $n\times n$ matrix $D$, \emph{seriation in the presence of 
errors} is to find a Robinsonian matrix $R$ 
that minimizes the error defined as: $\max ||D_{ij}-R_{ij}||$ over all 
$i$ and $j$ in $\{1,2,3,\ldots,n\}$. Chepoi et al. in 
\cite{chepoi2009seriation} proved that seriation in the presence of 
errors is an NP-Hard problem. 
Later in \cite{chepoi2011seriation}, Chepoi and Seston gave a factor 16 
approximation algorithm. Fortin in \cite{fortin2017robinsonian} surveyed
the challenges for Robinsonian matrix recognition. 

The SCFE problem was first introduced by Kermarrec and Thraves in
\cite{kermarrec2011can}. Besides the introduction of the SCFE problem, 
the authors of  \cite{kermarrec2011can}
also characterized the set of complete signed graphs with a valid 
drawing
in $\mathbb{R}$ and presented a polynomial time recognition algorithm. 
Later, Cygan et al. in \cite{cygan2015sitting}
proved that the SCFE problem is NP-Complete if it is not restricted to 
complete signed graphs.  Moreover, they gave a different 
characterization
of the complete signed 
graphs with a valid drawing in $\mathbb{R}$. Actually, the authors of 
\cite{cygan2015sitting} proved that a complete signed graph has a valid 
drawing in $\mathbb{R}$ if and only if its positive subgraph is a unit 
interval graph.
The SCFE problem in the real line also was studied as an optimization 
problem by Pardo et al. in \cite{pardo2015embedding}. In that work, the 
authors defined 
as an error a violation of the inequality in Definition 
\ref{def:validdrawing} and provided optimization algorithms that 
construct a drawing attempting to minimize the number of errors.

The SCFE problem also has been studied for different metric spaces.
First, Benitez et al. in \cite{benitez2018sitting} studied the 
SCFE
problem in the circumference. 
The authors of that work proved that the SCFE problem in the 
circumference is NP-Complete and gave a characterization of the complete
signed graphs with a valid drawing. 
Indeed, they showed that a complete signed graph has a valid drawing in 
the circumference if and only if its positive subgraph is a proper 
circular arc graph.  
Later, Becerra in \cite{becerra2018trees} studied the SCFE problem in 
trees. The main result of her work was to prove that a complete signed 
graph $G$ has a valid drawing in a tree if and only if its  positive 
subgraph is strongly chordal. 

Spaen et al. in \cite{spaen2017rk} studied the SCFE problem from a 
different perspective. They studied the problem of finding
$L(n)$, the smallest dimension $k$ such that any signed graph on $n$ 
vertices has a valid drawing in $\mathbb{R}^k$, with respect
to the Euclidean distance. They showed that $\log_5 (n-3) \leq L(n) \leq
n - 2$.

\paragraph{Our Contributions}
Our first contribution is to show that the Seriation and the SCFE 
problems are not the same. Indeed, we show that the SCFE problem implies a stronger condition 
than the Seriation problem. 
 In Lemma \ref{lem:validthenrobin}, 
we show that if a weighted graph $G=(V,E)$ has a valid drawing in $\mathbb{R}$, there is
a Robinson ordering of $V$. 
Nevertheless, in Lemma \ref{lem:VDmorethanSeriation}, 
we show that there is a weighted graph $G=(V,E)$ with a Robinson ordering
of $V$, but $G$ does not have a valid drawing in $\mathbb{R}$.

On the other hand, for each weighted graph $G$ with a Robinson ordering of its vertex set,
we construct a polyhedron defined by an inequality system 
$M(G)\mathbf{x}\leq \mathbf{b}$ which we use to provide a
characterization of the set of weighted graphs with a valid 
drawing in $\mathbb{R}$. 
Indeed, we show in Theorem \ref{thm:characterization} that a weighted 
graph $G$ has a valid drawing in $\mathbb{R}$ 
if and only if its polyhedron defined by $M(G)\mathbf{x}\leq \mathbf{b}$ is not 
empty.

Our first result applied to complete weighted graphs allows us to 
conclude in Corollary \ref{coro:polytimecomplete} that given a complete 
weighted graph $G$, determining whether $G$ has 
a valid drawing in $\mathbb{R}$, and finding one if applicable, can be 
done in polynomial time.

On the other hand, when the weighted graph is not complete, the previous
result does not apply anymore.
In Corollary 
\ref{coro:ReconIncompleteRob}, we state that
recognition of incomplete 
Robinsonian 
matrices is NP-complete. Furthermore,
using results shown in \cite{cygan2015sitting} by Cygan et al.,
we conclude, under the assumption of the Exponential Time Hypothesis,  the nonexistence of a subexponential-time algorithm
that determines if an incomplete similarity matrix is 
Robinsonian. 
Therefore, the construction of the polyhedron 
cannot be done in polynomial 
time (unless P=NP).  
Nevertheless, in Section \ref{sec:incompletegraphs},
we provide a recognition algorithm of $n\times n$ incomplete Robinsonian matrices with time complexity $O(n\cdot 2^{2n})$.
%

\section{Robinson Orderings and Valid Distance Drawings}\label{sec:ordersandDrawings} 
In this section we connect Robinson orderings and valid distance drawings. 
We start this section showing 
that, 
if a weighted graph has a valid drawing in 
$\mathbb{R}$ its  similarity matrix is Robinsonian. Therefore, 
having a Robinsonian similarity matrix is a necessary condition 
to have a valid distance drawing in $\mathbb{R}$.

\begin{lemma}\label{lem:validthenrobin}
Let $G=(V,E)$ be a weighted graph. If $G$ has a valid distance drawing in 
$\mathbb{R}$, then $V$ has a Robinson ordering.  
\end{lemma} 
\begin{proof}
Let $G=(V,E)$ be a  weighted graph with weight function $w$. Let 
$D:V\rightarrow \mathbb{R}$ be a valid distance drawing of $G$ in $\mathbb{R}$. 
The valid drawing $D$ determines an ordering on the set of vertices $V$.
Indeed, for $u$ and $v$ in $V$, we say that $u<_D v$ if $D(u) < D(v)$. 
We show that $A(G)$ is Robinson when it is written using the ordering determined by $D$ 
for its rows and columns.

Enumerate $V$ according to the ordering determined by $D$. 
Consider any integers $i,j, k,l$ such that 
$1\leq i < l \leq n$, $j,k \in [i,l]$,
and  $A(G)_{il} \neq *$, $A(G)_{ij} \neq *$ and 
$A(G)_{kl} \neq *$. First, we point our that 
$A(G)_{il} \leq 
A(G)_{ij}$, since $D$ is 
valid distance, and $d(D(i),D(l)) > d(D(i),D(j))$. Equivalently, since $D$ is valid distance and 
$d(D(i),D(l)) > d(D(k),D(l))$, we have $A(G)_{il} \leq A(G)_{kl}$. 
Therefore, $A(G)_{il} \leq \min \{A(G)_{ij},A(G)_{kl}\}$.
In conclusion, $A(G)$, the similarity matrix of $G$ written according to the ordering determined 
by $D$, is Robinson. 
Hence, the ordering determined by $D$ is Robinson. 

\end{proof}

Nevertheless, having a Robinson similarity matrix is not enough to have
a valid distance drawing in $\mathbb{R}$. 
In the next, lemma we show the existence of a weighted graph with a
Robinson ordering of its vertices but without valid distance drawing in 
$\mathbb{R}$. 
\begin{lemma}\label{lem:VDmorethanSeriation}
There exists a complete weighted graph $G$ with Robinson similarity matrix 
but without a valid distance drawing in $\mathbb{R}$.
\end{lemma}
\begin{proof}
Let $G$ be the complete weighted graph with vertex set $\{a,b,c,d,e\}$ 
and similarity matrix
\[
A(G)=
\begin{bmatrix}
5&2&2&1&1\\
2&5&3&2&1\\
2&3&5&4&1\\
1&2&4&5&5\\
1&1&1&5&5\\
\end{bmatrix}
\]
written with rows and columns ordered as $a,b,c,d,e$.
$A(G)$ is Robinson, nevertheless, we will show by contradiction 
that $G$ does not have a valid drawing in $\mathbb{R}$. 

Assume that $G$ has a valid drawing $D$ in $\mathbb{R}$. Since the order
$a,b,c,d,e$ of the rows and columns of $A(G)$ is the only one that 
makes $A(G)$ Robinson, $D$ has to be such that 
\begin{equation} \label{ineqorder}
D(a) < D(b) < D(c) < D(d) < D(e).
\end{equation}

Since $D$ is a valid drawing, the following 
inequalities hold: 
\begin{eqnarray}
D(b) - D(a) &>& D(c) -D(b) \label{ineq1}\\
D(e)-D(b) &>& D(b) - D(a) \label{ineq2}\\
D(c)-D(b) &>& D(d) - D(c) \label{ineq3}\\
D(e)-D(c) &>& D(c)-D(a) \label{ineq4}\\
D(d)-D(c) &>& D(e)-D(d) \label{ineq5}.
\end{eqnarray}

Without loss of generality, assume that $D(a)=0$. 
Then, from inequalities \eqref{ineqorder} and \eqref{ineq1} we obtain:
\begin{equation}\label{ineqcontradic1}
D(b) < D(c) <2D(b).    
\end{equation}
On the other hand, from inequalities \eqref{ineq4} and \eqref{ineq5},
we obtain $2D(c) < D(e) < 2D(d)-D(c)$, which implies:
\begin{equation}\label{ineqpivot}
    3D(c) < 2D(d).
\end{equation}
Finally, inequality \eqref{ineq3} is equivalent to $2D(d) < 
4D(c)-2D(b)$, 
which, together with \eqref{ineqpivot}, implies $2D(b) < 
D(c)$. 
But, the last inequality contradicts inequality \eqref{ineqcontradic1}.
\end{proof}

\section{The Weighted SCFE Problem in the line}\label{sec:realline}
The goal of this section is to find a solution for the weighted SCFE 
problem in the real line.  
We 
transform
the weighted SCFE problem in the real line 
into the problem of finding a point in a convex polyhedron.
Actually, given a weighted graph $G$, we construct a convex 
polyhedron defined by an inequality system $M(G)\mathbf{x}\leq \mathbf{b}$, where 
each point $\mathbf{x}=(x_1,x_2,\ldots,x_n)$ in the convex polyhedron is
a valid drawing of $G$ in $\mathbb{R}$. Indeed, for any given 
$\mathbf{x}$ in the polyhedron, each variable $x_i$ 
represents the position of vertex $i$ in $\mathbb{R}$ for that valid 
drawing. 
Therefore, finding a point in the polyhedron is 
equivalent to find a valid drawing for $G$ in $\mathbb{R}$.


We first remark that if a given weighted graph $G$ has a valid drawing 
in
$\mathbb{R}$, it actually has an infinite number of them. 
Indeed, given a valid drawing in $\mathbb{R}$ for a weighted graph $G$, 
one can obtain a different valid drawing for the same graph by summing 
or
multiplying 
each vertex position by any positive constant. The second case (when 
each
position is multiplied by a positive constant) is important for us, 
because 
it allows us to state the following lemma. 
\begin{lemma}\label{lem:foranyepsilon}
Let $G=(V,E)$ be a weighted graph with a valid distance drawing in $\mathbb{R}$. Then, 
for any $\epsilon > 0$ there exists a valid distance drawing $D_\epsilon$ of $G$ 
in 
$\mathbb{R}$ such 
that: 
 \[
 \min_{u,v \in V } |D_\epsilon(u)-D_\epsilon(v)| \geq \epsilon.
 \]
\end{lemma}
\begin{proof}
Let $G$ be a weighted graph with a valid drawing $D$ in $\mathbb{R}$. We
consider without loss of generality that $v$ is labeled according to the ordering determined by $D$, i.e., $1<_D2<_D3<_D\ldots <_D n$. 
Consider any $\epsilon > 0$. Let $\delta = \min_{ 1 \leq i < n } 
D(i+1)-D(i)$ be the minimum distance between 
two consecutive vertices in the drawing. Multiply every $D(i)$ 
by $\epsilon/\delta$. Therefore, we obtain a new valid drawing 
$D_\epsilon$ 
defined as $D_\epsilon(i)=\epsilon D(i)/\delta$, 
such that $\min_{u,v \in V } |D_\epsilon(u)-D_\epsilon(v)|\geq \epsilon$.

\end{proof}

Now, we proceed with the construction of the matrix $M(G)$ and the 
vector $\mathbf{b}$ of the inequality system $M(G)\mathbf{x}\leq 
\mathbf{b}$ that defines our polyhedron.
By Lemma \ref{lem:validthenrobin}, the ordering of the vertex set of $G$ defined by a valid 
drawing is Robinson. 
Hence, we may assume that $A(G)$ is Robinson.  

If we want to construct a valid drawing $D$ in $\mathbb{R}$ 
for $G$, the vertices should be ordered in the same way as the rows and 
columns of $A(G)$. 
Hence, if the $i$-\emph{th} row (or column) of $A(G)$ contains the 
similarities of vertex $i \in \{1,2,3,\ldots, n\}$, then $D$ has to be so that
$D(1) < D(2) < \cdots < D(n)$. Therefore, we want $x_1<x_2<\cdots 
<x_n$. Now, considering Lemma \ref{lem:foranyepsilon}, we write the 
following set of restrictions for any $\epsilon > 0$:
\begin{equation}\label{restriction:order}
x_i-x_{i+1} \leq -\epsilon, \quad \forall i \in \{1,2,3,\ldots,n-1\}.
\end{equation}
These restrictions are called \emph{ordering restrictions}. 
Lemma \ref{lem:foranyepsilon} allows us to pick any $\epsilon > 0$ 
for these restrictions. We set a value for $\epsilon$, and use that value for the inequalities we define below. 

On the other hand, each row of $A(G)$ provides two types of restrictions.
We call these restrictions \emph{right with respect to left} and  
\emph{left with respect to right} restrictions. 
Right with respect to left restrictions are obtained as follows. For 
each
row $j$ and for every index $k > j$, let $i(k)$ be the largest index 
such
that $i(k) < j$ and $A(G)_{ji(k)} < A(G)_{jk}$.
Therefore, since $A(G)_{ji(k)} < A(G)_{jk}$, vertices $j$ and $k$ are
more similar between them than vertices $j$ and $i(k)$. Hence, in any 
valid drawing  $D$ it must occur $D(k)-D(j) < D(j)-D(i(k))$.
We transform this strict inequality into the following restriction for a 
sufficiently small $\epsilon > 0$: 
\begin{equation}\label{restriction:rightwrtleft}
x_{i(k)}-2x_j+x_k  \leq -\epsilon, \quad \forall j \in 
\{2,3,\ldots,n-1\}
\mbox{ and } \forall  k>j.
\end{equation}

Left with respect to right restrictions are symmetrical to the previous 
restriction. For each row $j$ and for every index $i < j$, let $k(i)$ be
the smallest index such that $j < k(i)$ and $A(G)_{ji} >A(G)_{jk(i)}$.
Therefore, since $A(G)_{ji} >A(G)_{jk(i)}$, vertices $i$ and $j$ are 
more
similar between them than vertices $j$ and $k(i)$. Hence, in any valid 
drawing  $D$, it must occur $D(j)-D(i) < D(k(i))- D(j)$.
We transform this strict inequality into the following restriction for a 
sufficiently small $\epsilon > 0$: 
\begin{equation}\label{restriction:leftwrtright}
-x_{i}+2x_j-x_{k(i)}  \leq -\epsilon, \quad \forall j \in 
\{2,3,\ldots,n-1\} \mbox{ and } \forall  i<j.
\end{equation}

It is worth mentioning that some of the inequalities described in 
equations (\ref{restriction:rightwrtleft}) and 
(\ref{restriction:leftwrtright}) may be obtained from inequalities 
presented in Equation (\ref{restriction:order}) and different 
inequalities described in equations (\ref{restriction:rightwrtleft}) and
(\ref{restriction:leftwrtright}). 
Hence, some restrictions may be redundant. In an 
attempt to keep the presentation of this document 
clean and simple, we omit a discussion
in this regard. It is worth mentioning though that it does not impact
the
results of this document.
%

\begin{figure}[t]
\begin{center}
\resizebox{0.9\textwidth}{!}{\includegraphics*{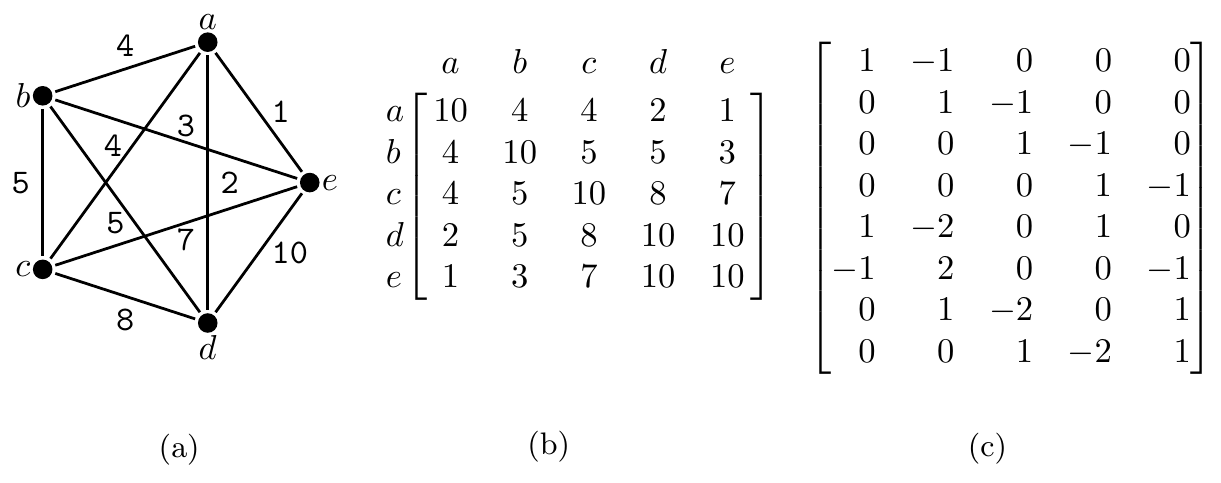}}
\caption{Example of a complete weighted graph, its similarity matrix, 
and
its corresponding matrix of restrictions. Subfigure (a) shows a complete
weighted graph. Subfigure (b) shows its Robinson similarity matrix. 
It also shows the order of the vertices in which the 
similarity matrix is written. Subfigure (c) shows the restriction matrix
for the weighted graph in Subfigure (a). In the first 4 rows appear the 
ordering restrictions. Rows five and six show the right with respect to 
left and left with respect to right restrictions for vertex $b$. Rows 
seven and eight show right with respect to left restrictions for
vertices
$c$ and $d$, respectively.}
\label{fig:example1}
\end{center}
\end{figure}

Given a weighted graph $G$ with $n$ vertices, the \emph{matrix
of restrictions} of $G$ (or the \emph{matrix of coefficients} of $G$),
denoted by  
$M(G)$, is the matrix that 
includes the $n-1$ ordering restrictions, the at most $(n-2)(n-1)/2$ 
right with respect to left restrictions, and  
the at most  $(n-2)(n-1)/2$ left with respect to right restrictions. 
In total, the matrix $M(G)$ has $h\leq (n-1)^2$ rows and $n$ columns.
On the other hand, the vector $\mathbf{b}$ is a $h \times 1$ vector 
with a $-\epsilon$ in every entry. 
An example of a weighted graph, its Robinson similarity matrix, 
and its corresponding matrix of restrictions is given in Figure 
\ref{fig:example1}.

Now, we show that for any weighted graph $G$ with Robinson similarity matrix, 
the convex polyhedron defined by $M(G)\mathbf{x}\leq \mathbf{b}$
is not empty if and only if $G$ has a valid drawing in $\mathbb{R}$.

\begin{theorem}\label{thm:characterization}
Let $G$ be a  weighted graph with Robinson similarity matrix $A(G)$. Let
$M(G)$ be the $h \times n$ matrix of restrictions of $G$ obtained from $A(G)$.
Let $\mathbf{b}$ be the 
$h \times 1$ vector 
with $-\epsilon < 0$ in every entry. Then, 
$G$ has a valid distance drawing in $\mathbb{R}$ if and only if
the polyhedron defined by $M(G)\mathbf{x} \leq \mathbf{b}$ 
is not empty. 
\end{theorem}
\begin{proof}
Let $G$ be a weighted graph with valid distance drawing in
$\mathbb{R}$. 
Let $D$ be a valid drawing of $G$ in $\mathbb{R}$.
Label the vertices of $G$ according to the order determined by $D$, i.e., 
the left most vertex in $D$ is vertex $1$, the next vertex is 
vertex $2$ and so on until vertex $n$. 
By construction of $M(G)\mathbf{x}\leq \mathbf{b}$,
for any $\epsilon > 0$, $D$ can be scaled to a valid drawing 
$D'$ such that
the vector $(D'(1),D'(2),\ldots,D'(n))$ belongs to the polyhedron 
$M(G)\mathbf{x}\leq \mathbf{b}$.

On the other hand, assume that the polyhedron 
$M(G)\mathbf{x}\leq \mathbf{b}$ 
is not empty. Let $x=(x_1,x_2\ldots,x_n)$ be a point in
$M(G)\mathbf{x}\leq \mathbf{b}$.
Label the vertices of $G$ according to the columns of its Robinson similarity matrix $A(G)$, 
i.e., vertex $i$ is the vertex 
corresponding to the $i$-\emph{th} column of $A(G)$.
Now, consider the drawing $D$ of $G$ in $\mathbb{R}$ defined as follows:
$D(i)=x_i$ for all $1\leq i \leq n$.

We show now that $D$ is valid distance. 
Assume that $D$ is not a valid distance drawing. Therefore, 
there exist three vertices $i,j$ and $k$ such that $A_{ij} < A_{ik}$,
but $|D(i)-D(j)| \leq |D(i)-D(k)|$. Note that the last inequality 
is not valid if 
$D(i) < D(k) < D(j)$ or if $D(j) < D(k) < D(i)$, therefore, these 
cases are discarded.
If $D(i)< D(j)<D(k)$ or
$D(k)<D(j)<D(i)$, there is a  
contradiction since $A_{ij}< A_{ik}$, and, in that case, $A(G)$
would not be Robinson. 

Assume that $D(j) < D(i) < D(k)$.
Therefore, $|D(i)-D(j)| \leq |D(i)-D(k)|$ becomes 
$D(i) -D(j) \leq D(k)-D(i)$, or equivalently, $0 \leq D(j)-2D(i)+D(k)$.
Nevertheless, since  $A_{ij} < A_{ik}$,  
the right with respect to left restriction $x_j - 2x_i +x_k \leq 
-\epsilon$ is included in 
$M(G)\mathbf{x} \leq \mathbf{b}$. Therefore, since $D$ comes 
from a point in $M(G)\mathbf{x} \leq \mathbf{b}$,
$D(j)-2D(i)+D(k) \leq -\epsilon$, which is a contradiction since 
$\epsilon > 0$. 

If we assume now  $D(k) < D(i) < D(j)$, then $|D(i)-D(j)| \leq 
|D(i)-D(k)|$ becomes $0 \leq -D(k)+2D(i)-D(j)$.
Nevertheless, since  $A_{ij} < A_{ik}$, the left with respect to 
right restriction $-x_k + 2x_i -x_j \leq 
-\epsilon$ is included in $M(G)\mathbf{x} \leq \mathbf{b}$. 
By equivalent arguments than 
before, we achieve a contradiction. 

Therefore, the condition $|D(i)-D(j)| \leq |D(i)-D(k)|$ 
is not possible, and hence, $D$ is a valid distance drawing. 
\end{proof}

If the valid drawings are restricted to be nonnegative, 
then the SCFE problem can be treated as a linear program. 
Because, if the polyhedron defined by the inequality system $M(G)\mathbf{x}\leq \mathbf{b}$ 
is not empty, 
there is always a point $\mathbf{x}$ in that polyhedron
with $x_1=0$.
Therefore, the SCFE problem is equivalent to find
$\min x_1$ subject to $M(G)\mathbf{x}\leq \mathbf{b}$, and 
nonnegative $\mathbf{x}$.

The last Theorem is stated for a weighted graph with a Robinson
similarity matrix. It is well known that a Robinsonian matrix may have many different Robinson orderings. Therefore, the reader may wonder which of these many Robinson orderings is the one that we can use to apply Theorem \ref{thm:characterization}. In the next part of this section, 
we answer that question for complete weighted graphs. Indeed, we show that any Robinson ordering will provide the same answer when Theorem \ref{thm:characterization} is applied.

Now we show that the existence of a valid
distance drawing in $\mathbb{R}$ 
corresponding to a Robinson ordering is consistent among all Robinson
orderings of the vertex set of a complete weighted graph.
In other words, we show that, 
given a complete weighted graph $G=(V,E)$, if there is a Robinson
ordering $\pi$ of $V$
and a valid distance drawing $D$ of $G$ such that the ordering induced
by $D$ is equal to $\pi$, 
then, for any Robinson ordering $\sigma$ of $V$ there exists a valid
distance drawing $\Sigma$
of $G$ such that the ordering induced by $\Sigma$ is equal to $\sigma$.

Given a Robinson ordering $\pi$ of $V$, we say that 
$\pi$ \emph{has a 
valid drawing} if
there is a valid distance drawing $D$ of $G=(V,E)$ in $\mathbb{R}$ such 
that 
the ordering induced by $D$ is equal to $\pi$.

To prove our result, we use the fact that all Robinson orderings for a complete similarity matrix
can be represented by a $PQ$-tree
(see \cite{prea2014optimal,booth1976testing}).
A $PQ$-tree on a set $V$ is a tree that represents a set of permutations of $V$. 
The nodes of a $PQ$-tree are of three types: leaves, that represent the elements
of $V$, $P$ nodes, and $Q$ nodes. The children of a $P$ node are not ordered, and any permutation 
of them is allowed. The children of a $Q$ node are ordered and that order can only be reversed.   
Hence, given a $PQ$-tree $\mathcal{T}$, we obtain a permutation of $V$ represented by 
$\mathcal{T}$ by applying one of the operations
allowed to each $P$ and $Q$ nodes, and then looking at the leaves to find the resultant permutation.
Let say $\pi$ is a permutation obtained in this way. If we modify the operation applied
to one node $\alpha$ of $\mathcal{T}$ and maintain all the other nodes equal, 
we say that the new permutation \emph{$\sigma$ is obtained from $\pi$ by modifying $\alpha$}. 

For a node $\alpha$ of a $PQ$-tree $\mathcal{T}$, 
we denote $\mathcal{T}(\alpha)$ the subtree of 
$\mathcal{T}$ with root $\alpha$ and by $S_\alpha$ the set of leaves of 
$\mathcal{T}(\alpha)$. 
A node of a $PQ$-tree is said to be \emph{basic} if all its children are 
leaves. Pr\'ea et al. in 
\cite{prea2014optimal} show that for every node $\alpha$, and for every
Robinson ordering $\pi$ 
represented in a $PQ$-tree $\mathcal{T}$, $S_\alpha$ is 
\emph{consecutive} according to $\pi$,
i.e., $S_\alpha=\{\pi_{l+1},\pi_{l+2},\ldots,\pi_{l+r}\}$ for some $0\leq l$ and  
$0\leq r \leq n-l$, where $\pi_i$ denotes the $i$-\emph{th} vertex according to $\pi$.
In addition, we denote by  $S_\alpha^I:=\pi_{l+1}$ and $S_\alpha^R:=\pi_{l+r}$  the first and the last elements of $S_\alpha$, respectively.
On the other hand, Pr\'ea et al. also show that 
\begin{equation}\label{eq:prea}
  w(\{u,x\})=w(\{v,x\}) \quad  \forall  u,v \in S_\alpha \mbox{ and } x 
\in V\setminus S_\alpha.
\end{equation}

We show now that operations on a $Q$-node of the $PQ$-tree
maintain the characteristic of having a valid distance drawing.
\begin{lemma}\label{lem:Qnodes}
Let $G=(V,E)$ be a complete weighted graph, $\mathcal{T}$ be the $PQ$-tree that 
represents all  Robinson orderings of $V$, and $\alpha$ be a 
$Q$-node 
of $\mathcal{T}$. 
Let $\pi$ and $\sigma$ be two Robinson orderings of $V$ such that
$\sigma$ is
obtained from $\pi$ by applying the operation associated to $\alpha$ (i.e., reversing the children of  $\alpha$). 
Then, $\pi$ has a valid distance drawing if and only if $\sigma$ has a valid distance drawing. 
\end{lemma}
\begin{proof}
Let $G$, $\mathcal{T}$, $\alpha$, $\pi$, and $\sigma$ be as in 
the statement of the lemma. 
Assume that $\pi$ has a valid drawing $D$. 
We show that $\sigma$ also 
has a valid drawing. 
It is worth noticing that this sense of the equivalence is
enough to show the lemma, since the opposite sense is shown by 
exchanging $\pi$ and $\sigma$ and repeating the analysis.  

Let $\pi_1,\pi_2,\ldots,\pi_n$ be the vertex set of $G$
ordered according to $\pi$.
The children of $\alpha$ can be $P$ nodes, $Q$ nodes or leaves.
Without loss of generality, we assume that all of them are $P$
nodes or $Q$ nodes, since a
leaf can be seen as a $Q$ node with a
single children that is the leaf. Let $\beta_1,\beta_2,\ldots,
\beta_h$ be $\alpha$'s children ordered according to $\pi$.
We use the following notation for the elements of  $S_{\beta_{t}}$: 
\[
\{\pi_{l+R_{t-1}+1},\pi_{l+R_{t-1}+2},\ldots,
\pi_{l+R_{t-1}+r_t}=\pi_{l+R_t}\},
\]
where $R_0=0$ and $R_{t} = R_{t-1}+ r_{t}$, and $|S_{\beta_{t}}|=r_t$.

Let $\tilde{\beta}_1,\tilde{\beta}_2,\ldots,
\tilde{\beta}_h$ be  $\alpha$'s children ordered according to $\sigma$.
We use the following notation for the elements of  $S_{\tilde{\beta}_{t}}$: 
\[
\{\sigma_{l+\tilde{R}_{t-1}+1},\sigma_{l+\tilde{R}_{t-1}+2},\ldots,\sigma_{l+\tilde{R}_{t-1}+\tilde{r}_t}=\sigma_{l+\tilde{R}_t}\},
\]
where $\tilde{R}_0=0$ and $\tilde{R}_{t} = \tilde{R}_{t-1} +\tilde{r}_t$ where $| S_{\tilde{\beta}_t}|=\tilde{r}_t =r_{h-t+1}$.

Note that $\tilde{R}_{t} = R_h-R_{h-t}$,  $R_{t} = \tilde{R}_h-\tilde{R}_{h-t}$ and $\tilde{R}_{h}=R_h$. In addition, $S_{\tilde{\beta}_t}= S_{\beta_{h-t+1}}$ and $S_{\beta_t}= S_{\tilde{\beta}_{h-t+1}}$.  In particular \[S_{\tilde{\beta}_t}^I=\pi_{l+R_{t-1}+1}= S_{\beta_{h-t+1}}^I= \sigma_{l+\tilde{R}_{h-t}+1}\] 
and \[S_{\tilde{\beta}_t}^R=\pi_{l+R_t}= S_{\beta_{h-t+1}}^R= \sigma_{l+\tilde{R}_{h-t+1}}.\]
Hence, 
\begin{eqnarray*}
 x \in S_{\tilde{\beta}_t }&\Longleftrightarrow& x \in \{\sigma_{l+\tilde{R}_{t-1}+1}, \ldots,\sigma_{l+\tilde{R}_t}\} 
 \\
 &\Longleftrightarrow& x \in S_{\beta_{h-t+1}} \\
 &\Longleftrightarrow& x \in \{\pi_{l+R_{h-t}+1},\ldots,\pi_{l+R_{t-h+1}} \}.
 \end{eqnarray*}
 
Assume that $\sigma_i=x$ is the $j$-element of $\tilde{\beta}_t$,  then 
\[
\sigma_i= \sigma_{l+\tilde{R}_{t-1}+j}  \Longrightarrow i=l+\tilde{R}_{t-1}+j \Longrightarrow j=i-l-\tilde{R}_{t-1}.
\]
On the other hand, $x$ is the $j$-element of  $\beta_{h-t+1}$, i.e., $x= \pi_{l+R_{h-t}+j}$. Hence: \[
\sigma_i=\pi_{l+R_{h-t}+j}= \pi_{l+R_{h-t}+i-l-\tilde{R}_{t-1}}= \pi_{i+R_{h-t}-R_h+R_{h-t+1}}.
\]
Therefore, $\sigma$ is as follows:
\[
\sigma_i=
\begin{cases}
\pi_i  &\mbox{if } i \in [1,l],\\
 \pi_{i+R_{h-t}-R_h+R_{h-t+1}} &\mbox{if } i \in [l+\tilde{R}_{t-1}+1,l+\tilde{R}_t], \ t \in \{1, \ldots, h\},\\
\pi_i  &\mbox{if } i \in [l+\tilde{R}_h+1,n].
\end{cases}
\]

We denote by $f:\{1, \ldots, n\}\to \{1, \ldots, n\}$ the  implicit bijection defined above in $\sigma_{i}=\pi_{f(i)}$. 
Next, we recursively define the  drawing $\Sigma$ such that $\Sigma(\sigma_{l+1})=D(\pi_{l+1})$ and $\forall i \in \{l+2, \ldots, n\}$: 
 \[
 d(\Sigma(\sigma_i), \Sigma(\sigma_{i-1}) )=
 \begin{cases} 
 d(D(\pi_{f(i)}), D(\pi_{f(i)-1})) &\mbox{if } \sigma_i \in S_{\tilde{\beta}_t}\setminus \{S_{\tilde{\beta}_t}^I \}, \\  &\mbox{and } t \in \{1, \ldots, h\} \\
 d(D(S^R_{\beta_{h-t+2}}), D(S^I_{\beta_{h-t+1}})) &\mbox{if } \sigma_i = S_{\tilde{\beta}_t}^I, \\ &\mbox{and } t \in \{2, \ldots, h\}. 
 \end{cases}
\]
Thus, we define the following drawing $\Sigma$ that induces $\sigma$: 
\begin{eqnarray*}
&&\Sigma(\sigma_i)= \\ 
&&\begin{cases}
D(\pi_i) &\mbox{if } i \in [1,l], \\
\Sigma(\sigma_{i-1}) + D(\pi_i)-D(\pi_{i-1})  &\mbox{if } i= l+1, \\
\Sigma(\sigma_{i-1}) + D(\pi_{f(i)}) - D(\pi_{f(i)-1}) & \mbox{if } i \in [l+\tilde{R}_{t-1}+2, l+\tilde{R}_t], \\ &\mbox{and } t \in [1, h],\\
\Sigma(\sigma_{i-1}) + D(\pi_{f(l+\tilde{R}_{t-2}+1)})-D(\pi_{f(l+\tilde{R}_t)}) & \mbox{if } i = l+\tilde{R}_{t-1}+1, \\ &\mbox{and } t \in [2, h],\\
D(\pi_i)& \mbox{if }  i \in [l+\tilde{R}_h+1,n].\\
\end{cases} 
\end{eqnarray*}

Observe that: 
\begin{enumerate}
\item $\Sigma(\sigma_{l+1})=\Sigma(\sigma_{l}) + D(\pi_{l+1}) - D(\pi_{l})= D(\pi_{l})+ D(\pi_{l+1}) - D(\pi_{l})=D(\pi_{l+1}),$
\item  $i = l+\tilde{R}_{t-1}+1 \Longleftrightarrow \sigma_i= S_{\tilde{\beta}_t}^I$, $\pi_{f(l+\tilde{R}_{t-2}+1)}=S_{\tilde{\beta}_{t-1}}^I $ and $\pi_{f(l+\tilde{R}_t)}=S_{\tilde{\beta}_t}^R$.
\end{enumerate}

We show now that $\Sigma$ is a valid distance drawing
for $G$ that induces $\sigma$ as an ordering. By contradiction, assume that there exists a triplet $\sigma_i$, 
$\sigma_j$,
$\sigma_k$ that breaks Definition \ref{def:validdrawing}, 
for some values $1\leq i<j<k\leq n$.
If the three elements of the triplet do not belong to $S_\alpha$,  
their positions according to $\Sigma$ do not change. Therefore, the 
distances between them do not change. Hence, they cannot violate 
Definition \ref{def:validdrawing}. 
On the other hand, if the three elements of the triplet belong to $S_{\beta_t}$ for some $1\leq t \leq h$, 
as we have seen in the previous paragraph, they maintain their distances. Therefore, 
they cannot violate 
Definition \ref{def:validdrawing}. 

Assume that only one element of the triplet belongs to $S_\alpha$. For instance, $\sigma_i\in S_\alpha$, while $\sigma_j\notin S_\alpha$ and $\sigma_k \notin S_\alpha$.
On the other hand, assume that 
$d(\Sigma(\sigma_i),\Sigma(\sigma_j)) \geq d(\Sigma(\sigma_j), \Sigma(\sigma_k))$, while $w(\sigma_i,\sigma_j) > w(\sigma_j,\sigma_k)$. Since $\sigma_i$ and $\sigma_{l+1}$
belong to $S_\alpha$ and $\sigma_j\notin S_\alpha$, 
Equation \eqref{eq:prea} implies
$w(\sigma_{l+1},\sigma_j) = w(\sigma_i,\sigma_j)$. Since, 
$d(\Sigma(\sigma_{l+1}),\Sigma(\sigma_j))\geq d(\Sigma(\sigma_i),\Sigma(\sigma_j))$, the triplet 
$\sigma_{l+1},\sigma_j,\sigma_k$ also breaks Definition 
\ref{def:validdrawing}.

Now,
\begin{eqnarray*}
d(\Sigma(\sigma_{l+1}),\Sigma(\sigma_j)) = \Sigma(\sigma_j) - \Sigma(\sigma_{l+1}) &&= D(\pi_j) - D(\pi_{l+1})\\ &&= d(D(\pi_{l+1}),D(\pi_j)),
\end{eqnarray*}
and
\[
d(\Sigma(\sigma_j), \Sigma(\sigma_k)) = \Sigma(\sigma_k) - \Sigma(\sigma_j) = D(\pi_k)- D(\pi_j) = d(D(\pi_j),D(\pi_k)).
\]
Therefore, $d(D(\pi_{l+1}),D(\pi_j))\geq d(D(\pi_j),D(\pi_k))$. 
On the other hand, \[w(\sigma_{l+1},\sigma_j) = w(\pi_{l+R_{h-1}+1},\pi_j),\] since $\sigma_l=\pi_{l+R_{h-1}+1}$ and $\sigma_j=\pi_j$. But, 
Equation \eqref{eq:prea} implies
\[w(\pi_{l+R_{h-1}+1},\pi_j) = w(\pi_{l+1},\pi_j),\] since $\pi_{l+R_{h-1}+1}$ and $\pi_l$ belong to $S_\alpha$, and $\pi_j \notin S_\alpha$. 
In conclusion, \[w(\pi_{l+1},\pi_j) = w(\pi_{l+R_{h-1}+1},\pi_j) = w(\sigma_l,\sigma_j)  > w(\sigma_j,\sigma_k) = w(\pi_j,\pi_k).\] 
Hence, the triplet $\pi_{l+1},\pi_i,\pi_k$ breaks Definition 
\ref{def:validdrawing} in $D$, which is a contradiction since $D$ is valid distance. 
The analysis of all cases when one element of 
the triplet belongs to $S_\alpha$ are equivalent to this analysis 
by showing that either $\sigma_{l+1},\sigma_j,\sigma_k$ or  $\sigma_{l+R_h},\sigma_j,\sigma_k$ also breaks Definition 
\ref{def:validdrawing}. We omit them to simplify the presentation. 

Assume now that two elements of the triplet belong to $S_\alpha$. Say $\sigma_i$ and $\sigma_j$ belong to $S_\alpha$, while $\sigma_k$ does not belong 
to $S_\alpha$. Assume as well 
that \[d(\Sigma(\sigma_i),\Sigma(\sigma_j)) \geq d(\Sigma(\sigma_j), \Sigma(\sigma_k)),\] while $w(\sigma_i,\sigma_j) > w(\sigma_j,\sigma_k)$.

Since $\sigma$ is Robinson, we have that:
\[
w(\sigma_{l+1},\sigma_{l+R_h}) \geq w(\sigma_{l+1},\sigma_j) \geq w(\sigma_i,\sigma_j),
\]
and
\[
w(\sigma_j,\sigma_k) \geq w(\sigma_{l+R_h},\sigma_k).
\]
Therefore, $w(\sigma_{l+1},\sigma_{l+R_h}) > w(\sigma_{l+R_h},\sigma_k)$.
On the other hand, \[d(\Sigma(\sigma_{l+1}),\Sigma(\sigma_{l+R_h})) \geq d(\Sigma(\sigma_{l+1}),\Sigma(\sigma_j)) \geq d(\Sigma(\sigma_i),\Sigma(\sigma_j))\] and  
\[ d(\Sigma(\sigma_j), \Sigma(\sigma_k)) \geq d(\Sigma(\sigma_{l+R_h}), \Sigma(\sigma_k)).
\]

Hence, \[d(\Sigma(\sigma_{l+1}),\Sigma(\sigma_{l+R_h})) \geq d(\Sigma(\sigma_{l+R_h}), \Sigma(\sigma_k))\] and 
the triplet $\sigma_{l+1},\sigma_{l+R_h},\sigma_k$ also breaks Definition \ref{def:validdrawing}.
Now,
\begin{eqnarray*}
d(\Sigma(\sigma_{l+1}),\Sigma(\sigma_{l+R_h})) &=& 
\Sigma(\sigma_{l+R_h})-\Sigma(\sigma_{l+1})\\
&=& D(\pi_{l+R_h})-D(\pi_{l+1})\\ &=& d(D(\pi_{l+1}),D(\pi_{l+R_h})),
\end{eqnarray*}
and
\begin{eqnarray*}
d(\Sigma(\sigma_{l+R_h}), \Sigma(\sigma_k)) &=& \Sigma(\sigma_k) - \Sigma(\sigma_{l+R_h})\\ &=& D(\pi_k) - D(\pi_{l+R_h})\\ &=& d(D(\pi_{l+R_h}),D(\pi_k)).
\end{eqnarray*}
Therefore, \[d(D(\pi_{l+1}),D(\pi_{l+R_h})) 
\geq d(D(\pi_{l+R_h}),D(\pi_k).\] 
On the other hand, applying repeatedly 
Equation \eqref{eq:prea} 
we obtain:
\begin{eqnarray*}
w(\pi_{l+R_h},\pi_{l+1})&=&w(\sigma_{l+1},\sigma_{l+R_h})\\ &>& w(\sigma_{l+R_h},\sigma_k)\\ &=& w(\pi_{l+1},\pi_k)\\ &=& w(\pi_{l+R_h},\pi_k).
\end{eqnarray*}
Therefore, the triplet $\pi_{l+1},\pi_{l+R_h},\pi_k$
breaks Definition \ref{def:validdrawing}. Which is a contradiction
because $D$ is valid distance. The analysis of all cases when two elements of 
the triplet belong to $S_\alpha$ are equivalent to this analysis. We omit them to simplify the presentation. 

Finally, assume that the three elements 
of the triplet belong to $S_\alpha$.  
These three vertices cannot all be in 
$S_{\tilde{\beta}_t}$ for any $t$ since, in this 
case, the distances between them according to $\Sigma$ 
remain the same as in $D$. Therefore, they could not break 
Definition \eqref{def:validdrawing}.
Assume that $\sigma_i$ and $\sigma_j$ belong to
$S_{\tilde{\beta}_t}$
and $\sigma_k$ belongs to $S_{\tilde{\beta}_p}$ with 
$1\leq p < t \leq h$.
Assume as well that 
$d(\Sigma(\sigma_i),\Sigma(\sigma_j)) \geq d(\Sigma(\sigma_j), 
\Sigma(\sigma_k))$, while $w(\sigma_i,\sigma_j) > 
w(\sigma_j,\sigma_k)$.

It is worth noticing, that $\sigma$ inverts the ordering of 
the children of $\alpha$. 
Indeed, using the relationship  $\tilde{\beta}_t= \beta_{h-t+1}$, 
we have: 
\begin{eqnarray*}
S_{\tilde{\beta}_t} &=& \{\sigma_{l+\tilde{R}_{t-1}+1},\ldots,\sigma_{l+\tilde{R}_t}\}= S_{\beta_{h-t+1}}  = \{\pi_{l+R_{h-t}+1},\ldots,\pi_{l+R_{h-t+1}}\}
\end{eqnarray*}
and
\begin{eqnarray*}
S_{\tilde{\beta}_p} &=& \{\sigma_{l+\tilde{R}_{p-1}+1},\ldots,\sigma_{l+\tilde{R}_p}\}= S_{\beta_{h-p+1}}  = \{\pi_{l+R_{h-p}+1},\ldots,\pi_{l+R_{h-p+1}}\}
\end{eqnarray*}

By Equation \eqref{eq:prea} \[w(\sigma_i,\sigma_j)>w(\sigma_j,\sigma_k)=w(\sigma_k,\sigma_{l+\tilde{R}_{t-1}+1})=w(\sigma_{l+\tilde{R}_{t-1}+1},\sigma_{l+\tilde{R}_p}).\]
 Since $D$ is a valid drawing, \[d(D(\pi_{i}),D(\pi_{j}))<d(D(\pi_{l+R_{h-t}+1}),D(\pi_{l+R_{h-p+1}})).\]
Thus, 
\begin{eqnarray*}
d(\Sigma(\sigma_{i}),\Sigma(\sigma_{j}))&=d(D(\pi_{i}),D(\pi_{j}))<d(D(\pi_{l+R_{h-t}+1}),D(\pi_{l+R_{h-p+1}})) \\
&= d(\Sigma(\sigma_{l+\tilde{R}_t}), \Sigma(\sigma_{l+\tilde{R}_{p-1}+1})) \leq d(\Sigma(\sigma_j), 
\Sigma(\sigma_k)),
\end{eqnarray*}
 which is a contradiction. 
The analysis of all cases when the three elements of 
the triplet belong to $S_\alpha$ are equivalent to this analysis. We omit them for simplicity of the presentation. 
\end{proof}

Now, we show that any operation on a $P$-node of the $PQ$-tree
maintains the characteristic of having a valid distance drawing.
\begin{lemma}
Let $G=(V,E)$ be a complete weighted graph, $\mathcal{T}$ be the
$PQ$-tree that 
represents all  Robinson orderings of $V$, and $\alpha$ be a
$P$-node 
of $\mathcal{T}$. 
Let $\pi$ and $\sigma$ be two Robinson orderings of $V$ such
that 
$\sigma$ is
obtained from $\pi$ via a permutation of the children of $\alpha$.
Then, $\pi$ has a valid drawing if and only if $\sigma$ has a
valid 
drawing. 
\end{lemma}
\begin{proof}
We use Lemma \ref{lem:Qnodes} in this proof, since 
this lemma is a particular case of the previous lemma. 

Let $G$, $\mathcal{T}$, $\alpha$, $\pi$, and $\sigma$ be as in 
the statement of the lemma. 
Assume that $\pi$ has a valid drawing $D$. 
We will show that $\sigma$ also 
has a valid drawing. 
It is worth noticing that, as in the previous lemma, this sense of the equivalence is enough to show the lemma. 

The children of $\alpha$ can be $P$ nodes, $Q$ nodes or leaves.
Without loss of generality, we assume that all of them are $P$
nodes or $Q$ nodes, since a
leaf can be seen as a $Q$ node with a
single children that is the leaf. Let $\beta_1,\beta_2,\ldots,
\beta_h$ be $\alpha$'s children ordered according to $\pi$. 
It is enough for us to consider that $\sigma$ is obtained via 
a transposition of two 
consecutive children of $\alpha$, since any permutation of
$S_\alpha$ 
can be 
expressed as a product of
transpositions of consecutive elements in $S_\alpha$.
Let $\beta_t$ and $\beta_{t+1}$ be the two children to be 
transposed that produce $\sigma$.
Now, let $\mathcal{T'}$ be a $PQ$-tree obtained from $\mathcal{T}$ by replacing $\beta_t$ and $\beta_{t+1}$
with a $Q$ node that has as children $\beta_t$ and $\beta_{t+1}$.
Now, it is worth noticing that $\mathcal{T'}$
encodes $\pi$ and $\sigma$. Furthermore, Lemma \ref{lem:Qnodes}
implies that if $\pi$ has a valid drawing, then $\sigma$
also has a valid drawing, and the proof is completed. 
\end{proof}

With these two lemmas together we can state the result. 
\begin{theorem}\label{thm:one-all}
Let $G=(V,E)$ be a complete weighted graph and $\pi$ be any Robinson ordering of $V$. 
If $\pi$ has a valid drawing, all Robinson orderings of $V$ have a valid drawing.   
\end{theorem}

Since
complete 
Robinsonian matrices 
can be recognized in time $O(n^2)$, it is possible to 
construct
the matrix $M(G)$ in polynomial time when $G$ is complete. 
Therefore, we can state the 
following corollary. 
\begin{corollary}\label{coro:polytimecomplete}
Let $G$ be a complete weighted graph. Deciding whether $G$ has a valid 
drawing in $\mathbb{R}$ can be done in polynomial time. Moreover, a 
valid drawing for $G$ in $\mathbb{R}$ can be computed also in 
polynomial time if such drawing exists.
\end{corollary}

\section{The Weighted SCFE Problem for Incomplete 
Weighted Graphs}\label{sec:incompletegraphs}
The construction presented in the previous section can also be applied to 
incomplete weighted graphs. The only requirement is that the matrix
$A$ is presented as a Robinson matrix. 
Nevertheless, we will see in this section that, 
if the condition of being complete is not requested for the weighted 
graph, it is not possible to determine in polynomial 
time whether its similarity matrix is Robinsonian or not, unless P=NP.
Despite this bad result, we present an exponential-time algorithm 
to recognize incomplete Robinsonian matrices. Hence, once this recognition 
has been done, we can apply the tools developed in the previous section 
to solve the weighted SCFE problem.
It is worth noticing that, in the case of incomplete weighted graphs, we do not have a result equivalent to Theorem \ref{thm:one-all}. Hence, we cannot guarantee that after applying 
the methodology developed in the previous section to incomplete weighted graphs, we will obtain a definitive answer.

Cygan et al. in \cite{cygan2015sitting}
proved the NP-Completeness of the particular case of the SCFE 
problem where the weight in the edges can only take values
$+1$ or $-1$. On the other hand, Kermarrec and Thraves proved  
in \cite{kermarrec2011can} the following theorem rephrased in 
our own words.
 \begin{theorem}\label{thm:kerthr}[Lemmas 3 and 4 in  \cite{kermarrec2011can}]
 Let $G=(V,E)$ be a weighted graph and $w:E\rightarrow 
 \{+1,-1\}$ be its weight function. 
Then, $G$ has a valid distance drawing in $\mathbb{R}$ 
if and only if there exists an ordering $\pi$ of $V$ 
that verifies the following two conditions: 
\begin{enumerate}
\item For all $i<j<k$ such that $\{\pi_i,\pi_k\}$ and 
$\{\pi_j,\pi_k\}$ belong to $E$, 
\[w(\{\pi_j,\pi_k\})=-1 \quad \implies \quad 
w(\{\pi_i,\pi_k\})=-1.\]
\item For all $i<j<k$ such that $\{\pi_i,\pi_k\}$ and 
$\{\pi_i,\pi_j\}$ belong to $E$,
\[w(\{\pi_i,\pi_j\})=-1\quad \implies \quad w(\{\pi_i,\pi_k\})=-1.\]  
\end{enumerate}
\end{theorem}
In other words, a weighted graph $G=(V,E)$ with weights $+1$ or $-1$ 
has a valid distance drawing in the line 
if and only if there exists an ordering $\pi$ of $V$ such that $A^{\pi}(G)$ is 
Robinson.
Hence, using these two results, we can state the following corollary.
\begin{corollary}\label{coro:ReconIncompleteRob}
Let $G$ be an incomplete weighted graph and $A(G)$ be its similarity matrix. 
Deciding whether $A(G)$ is Robinsonian or not is a NP-Complete problem. 
\end{corollary}

Besides this negative result, Cygan et al. also proved in 
\cite{cygan2015sitting}  
the existence of a constant $C>0$  
such that no algorithm
solves the  SCFE problem in the line in time $O(2^{C(n+m)})$, unless the  
Exponential Time Hypothesis\footnote{The Exponential Time Hypothesis 
states that there exists a constant $C>0$ such that no algorithm solving the
3-CNF-SAT problem in $O(2^{CN})$ exists, where $N$ denotes the number of variables
in the input formula.} (ETH) fails.
Therefore, we conclude that, under the assumption of the ETH,
it is impossible to have a
subexponential-time algorithm 
that determines if the similarity matrix of an 
incomplete weighted graph is Robinsonian. 

On the positive side, 
using similar ideas
to those presented in \cite{cygan2015sitting}, 
we present an exponential-time algorithm 
that decides if a given incomplete similarity matrix 
is Robinsonian or not.
We start with the following definition.
\begin{definition}\label{def:good}
Let $A$ be an incomplete $n\times n$ similarity matrix. 
Let $\{V,U\}$ be a bipartition of the set $\{1,2,\ldots,n\}$.
We say that an element $u \in U$ is \emph{good} for $V$
if for all $k \in V$ and $p\in U$
\[
A_{ku} \neq * \, \land \,  A_{kp} \neq * \implies A_{ku}\geq A_{kp},
\]
and
\[
A_{pu} \neq *  \, \land \, A_{pk} \neq * \implies
A_{pu}\geq A_{pk}.
\]
\end{definition}

Now, we state the following result. 
\begin{lemma} 
Let $A$ be an incomplete $n\times n$ similarity matrix. 
$A$ is Robinsonian if and only if 
there exists an ordering $\pi=\pi_1, \pi_2,\ldots,\pi_n$ of
the set $\{1,2,\ldots,n\}$
such that for every 
$1\leq j \leq n-1$, the element $\pi_{j+1}$ is good for 
$\{\pi_1, \ldots,\pi_j \}$.
\end{lemma}
\begin{proof} 
Let $A$ be Robinsonian and  $\pi=\pi_1, \pi_2, \ldots, \pi_n$ be a Robinson
ordering of its rows and columns. Therefore, $A^\pi$ is Robinson. 
We show by contradiction that the ordering $\pi$ satisfies the conditions
of the lemma. Assume for instance that the first condition is 
broken for some $1\leq j \leq n-1$. 
Therefore, there exists  $1\leq k \leq j$ and $j+1 < p \leq n$ such that: 
$A_{\pi_k\pi_{j+1}} \neq *$, 
$A_{\pi_k\pi_p}\neq *$, and 
$A_{\pi_k\pi_{j+1}} < A_{\pi_k\pi_p}$. 
Since, 
$A_{\pi_k\pi_{j+1}}=A^\pi_{kj+1}$, $A_{\pi_k\pi_p}=A^\pi_{kp}$, and $j+1 < p$,
we obtain a contradiction with the fact that $A^\pi$ is Robinson.  
A similar conclusion can be drawn if we assume that 
the second condition of the lemma is 
broken for some $1\leq j \leq n-1$. 
Therefore, the necessary condition for $A$ to be  Robinsonian holds.
 
Assume now that there exists an ordering 
$\pi=\pi_1, \pi_2, \ldots,\pi_n$ such that
for every $1\leq j \leq n -1$, the element $\pi_{j+1}$ is good 
for the set $\{\pi_1,\ldots, \pi_ j\}$.
We show that $A^{\pi}$ is Robinson. 
Consider integers $i$ and $l$ such that $1\leq i < l \leq n$. 
Consider now any integer $j$ in $[i,l]$. 
Since the element $\pi_j$ is good for the set $\{\pi_1,\pi_2,\ldots,\pi_{j-1}\}$, 
the first condition of Definition \ref{def:good} implies that 
if $A^\pi_{ij}\neq *$ and $A^\pi_{il}\neq *$, then 
$A^\pi_{ij} \geq A^\pi_{il}$. Equivalently, 
since the element $\pi_k$ is good for the set $\{\pi_1,\pi_2,\ldots,\pi_{k-1}\}$,
if we consider an integer $k\in [i,l]$, the second condition of 
Definition \ref{def:good} implies that if 
$A^\pi_{kl}\neq *$ and $A^\pi_{il}\neq *$, then 
$A^\pi_{kl} \geq A^\pi_{il}$. 
Therefore, for any $1\leq i < l \leq n$ and $j,k \in [i,l]$, we have that 
$A^\pi_{il} \leq \min \{A^\pi_{ij},A^\pi_{kl}\}$. 
Hence, $A^\pi$ is Robinson. 
\end{proof}

Finally, to determine  if  an  incomplete similarity matrix $A$ is Robinsonian we
use the algorithm presented in \cite{cygan2015sitting} which 
proceeds as follows: construct a directed graph
$H=(S,F)$ where the vertex set $S$ is formed by  all  subsets of 
$\{1,2,\ldots, n\}$,  and for every $X,Y \subseteq \{1,2,\ldots,n\}$ there exists
an arc from $X$ to $Y$ if $Y\setminus X=\{i\}$ and $i$ is a good element for $X$.
The arc $(X,Y) \in F$ is labeled by $i$.  Thus, the existence of an ordering 
$\pi$ such that $A^{\pi}$ is Robinson 
is equivalent to the existence of a directed path from a vertex that is a
singleton of $\{1,2,\ldots, n\}$ to the vertex $\{1,2\ldots,n\}$. 
The existence of such directed path in $H$ can be
determined in time $O(n\cdot 2^{2n})$, using repeatedly the single source shortest path algorithm presented in \cite{THORUP2004330} for all the $n$ vertices representing a singleton.
Furthermore, the Robinson ordering $\pi$ is determined by the ordering 
of the labels along that directed path.

\section{Final Remarks}\label{sec:conclusions}
Interestingly, in this work we show that the Seriation and the
SCFE problems are not the same. Nevertheless, there are cases in which 
they are equivalent. For instance, an exhaustive analysis shows that
if a weighted graph has at most four vertices then its 
similarity matrix is Robinsonian if and only if it has a valid
drawing in $\mathbb{R}$. Whereas, in the proof of 
Lemma \ref{lem:VDmorethanSeriation}, we present a weighted graph
with five vertices where seriation is not sufficient. 

The Seriation and the SCFE problems are also equivalent if the
number of different weights is not too big. Indeed,
Theorem \ref{thm:kerthr} states that when the weight function can take only two values, seriation and the SCFE problem are equivalent.
Nevertheless, in the proof of 
Lemma \ref{lem:VDmorethanSeriation}, we exhibit an example of a weighted
graph  with five different weights where seriation is not 
enough. This final remark rises an interesting question, when this 
separation between the Seriation and the SCFE problems occurs?. Is the 
Seriation problem  
equivalent to the SCFE problem when the graph has four different 
weights?.

On the other hand, in Theorem \ref{thm:one-all} we have established 
that, when the weighted graph is complete, if one Robinson ordering 
has a valid distance drawing, then all Robinson orderings have one. 
Such a result is crucial to show that the SCFE problem has a polynomial 
time algorithm when the input is a complete weighted graph.
To prove this Theorem, we use the fact that there is a $PQ$-tree that 
encodes all Robinson orderings for a complete weighted graph. 
We do not have a result like 
that  for incomplete weighted graphs. Hence, it is not clear 
the existence of a result for incomplete weighted graphs equivalent to
Theorem \ref{thm:one-all}. We believe that this is a really interesting 
problem that remains open. Indeed, is it possible to encode all Robinson
Orderings for an incomplete weighted graph in a $PQ$-tree?.


\begin{thebibliography}{10}

\bibitem{becerra2018trees}
{\sc R.~Becerra and C.~Thraves~Caro}, {\em On the sitting closer to friends
  than enemies problem in trees and an intersection model for strongly chordal
  graphs}, arXiv preprint arXiv:1911.11494,  (2019).

\bibitem{benitez2018sitting}
{\sc F.~Ben{\'\i}tez, J.~Aracena, and C.~Thraves~Caro}, {\em The sitting closer
  to friends than enemies problem in the circumference}, arXiv preprint
  arXiv:1811.02699,  (2018).

\bibitem{booth1976testing}
{\sc K.~S. Booth and G.~S. Lueker}, {\em Testing for the consecutive ones
  property, interval graphs, and graph planarity using pq-tree algorithms},
  Journal of computer and system sciences, 13 (1976), pp.~335--379.

\bibitem{brusco2006branch}
{\sc M.~J. Brusco and S.~Stahl}, {\em Branch-and-{B}ound applications in
  combinatorial data analysis}, Springer Science \& Business Media, 2006.

\bibitem{chepoi1997recognition}
{\sc V.~Chepoi and B.~Fichet}, {\em Recognition of {R}obinsonian
  dissimilarities}, Journal of Classification, 14 (1997), pp.~311--325.

\bibitem{chepoi2009seriation}
{\sc V.~Chepoi, B.~Fichet, and M.~Seston}, {\em Seriation in the presence of
  errors: {NP}-hardness of $l_{\infty}$-fitting {R}obinson structures to
  dissimilarity matrices}, Journal of classification, 26 (2009), pp.~279--296.

\bibitem{chepoi2011seriation}
{\sc V.~Chepoi and M.~Seston}, {\em Seriation in the presence of errors: A
  factor 16 approximation algorithm for $l_{\infty}$-fitting {R}obinson
  structures to distances}, Algorithmica, 59 (2011), pp.~521--568.

\bibitem{cygan2015sitting}
{\sc M.~Cygan, M.~Pilipczuk, M.~Pilipczuk, and J.~O. Wojtaszczyk}, {\em Sitting
  closer to friends than enemies, revisited}, Theory of Computing Systems, 56
  (2015), pp.~394--405.

\bibitem{ding2004linearized}
{\sc C.~Ding and X.~He}, {\em Linearized cluster assignment via spectral
  ordering}, in Proceedings of the twenty-first international conference on
  Machine learning, ACM, 2004, p.~30.

\bibitem{fortin2017robinsonian}
{\sc D.~Fortin}, {\em Robinsonian matrices: {R}ecognition challenges}, Journal
  of Classification, 34 (2017), pp.~191--222.

\bibitem{hubert2001combinatorial}
{\sc L.~Hubert, P.~Arabie, and J.~Meulman}, {\em Combinatorial data analysis:
  Optimization by dynamic programming}, vol.~6, SIAM, 2001.

\bibitem{kermarrec2011can}
{\sc A.-M. Kermarrec and C.~Thraves~Caro}, {\em Can everybody sit closer to
  their friends than their enemies?}, in Proceedings of the 36th International
  Symposium on Mathematical Foundations of Computer Science, Springer, 2011,
  pp.~388--399.

\bibitem{laurent2017lex}
{\sc M.~Laurent and M.~Seminaroti}, {\em A {L}ex-{BFS}-based recognition
  algorithm for {R}obinsonian matrices}, Discrete Applied Mathematics, 222
  (2017), pp.~151--165.

\bibitem{laurent2017similarity}
{\sc M.~Laurent and M.~Seminaroti}, {\em Similarity-first search: {A} new
  algorithm with application to {R}obinsonian matrix recognition}, SIAM Journal
  on Discrete Mathematics, 31 (2017), pp.~1765--1800.

\bibitem{laurent2017structural}
{\sc M.~Laurent, M.~Seminaroti, and S.-i. Tanigawa}, {\em A structural
  characterization for certifying {R}obinsonian matrices}, The Electronic
  Journal of Combinatorics, 24(2) (2017).

\bibitem{liiv2010seriation}
{\sc I.~Liiv}, {\em Seriation and matrix reordering methods: An historical
  overview}, Statistical Analysis and Data Mining: The ASA Data Science
  Journal, 3 (2010), pp.~70--91.

\bibitem{graphsandgenes84}
{\sc B.~G. Mirkin and S.~N. Rodin}, {\em Graphs and {G}enes}, Springer-Verlag,
  1984.

\bibitem{pardo2015embedding}
{\sc E.~G. Pardo, M.~Soto, and C.~Thraves~Caro}, {\em Embedding signed graphs
  in the line}, Journal of Combinatorial Optimization, 29 (2015), pp.~451--471.

\bibitem{petrie1899sequences}
{\sc W.~M.~F. Petrie}, {\em Sequences in prehistoric remains}, Journal of the
  Anthropological Institute of Great Britain and Ireland,  (1899),
  pp.~295--301.

\bibitem{prea2014optimal}
{\sc P.~Pr{\'e}a and D.~Fortin}, {\em An optimal algorithm to recognize
  {R}obinsonian dissimilarities}, Journal of Classification, 31 (2014),
  pp.~351--385.

\bibitem{roberts69}
{\sc F.~S. Roberts}, {\em Indifference graphs}, in Proof Techniques in Graph
  Theory, F.~Harary, ed., Academic Press, New York, 1969, pp.~139--146.

\bibitem{robinson_1951}
{\sc W.~S. Robinson}, {\em A method for chronologically ordering archaeological
  deposits}, American Antiquity, 16 (1951), pp.~293 -- 301.

\bibitem{spaen2017rk}
{\sc Q.~Spaen, C.~Thraves~Caro, and M.~Velednitsky}, {\em The dimension of
  valid distance drawings of signed graphs}, Discrete \& Computational
  Geometry, 63 (2020), pp.~158--168.

\bibitem{THORUP2004330}
{\sc M.~Thorup}, {\em Integer priority queues with decrease key in constant
  time and the single source shortest paths problem}, Journal of Computer and
  System Sciences, 69 (2004), pp.~330 -- 353.

\bibitem{tien2008methods}
{\sc Y.-J. Tien, Y.-S. Lee, H.-M. Wu, and C.-H. Chen}, {\em Methods for
  simultaneously identifying coherent local clusters with smooth global
  patterns in gene expression profiles}, BMC bioinformatics, 9 (2008), p.~155.

\end{thebibliography}

\end{document}